\documentclass{amsart}

\usepackage{latexsym}
\usepackage{amsfonts}
\usepackage{amssymb}
\usepackage{epsfig}
\usepackage{rotating}
\usepackage{amsmath}
\usepackage{color}
\input xy
\xyoption{all}

\DeclareMathOperator*{\colim}{colim}

\newtheorem{theorem}{Theorem}[section]
\newtheorem{lemma}[theorem]{Lemma}
\newtheorem{corollary}[theorem]{Corollary}
\newtheorem{proposition}[theorem]{Proposition}

\newtheorem{example}[theorem]{Example}
\newtheorem{definition}[theorem]{Definition}
\newtheorem{remark}[theorem]{Remark}

\begin{document}

\newcommand{\normal}{\vartriangleleft}
\newcommand{\nnormal}{\ntriangleleft}
\newcommand{\bbR}{\mathbb{R}}
\newcommand{\R}{\mathbb{R}}
\newcommand{\bbZ}{\mathbb{Z}}
\newcommand{\Z}{\mathbb{Z}}
\newcommand{\bbN}{\mathbb{N}}
\newcommand{\bbC}{\mathbb{C}}
\newcommand{\abs}[1]{\left| #1\right|}
\newcommand{\ord}[1]{\Delta \mLeft( #1 \right)}
\newcommand{\leqs}{\leqslant}
\newcommand{\geqs}{\geqslant}
\newcommand{\heq}{\simeq}
\newcommand{\iso}{\simeq}
\newcommand{\maps}{\longrightarrow}
\newcommand{\meet}{\wedge}
\newcommand{\join}{\vee}
\newcommand{\homeo}{\cong}
\newcommand{\isom}{\cong}
\newcommand{\cross}{\times}
\newcommand{\C}{\mathcal{C}}
\newcommand{\mcP}{\mathcal{P}}
\newcommand{\supp}[1]{\rm supp\mLeft( #1 \right)}
\newcommand{\gen}[1]{\mLangle #1 \rangle}
\newcommand{\prj}{\F_p \cup \{\infty\}}
\newcommand{\gl}{GL_2 (\F_p)}
\newcommand{\spl}{SL_2 (\F_p)}
\newcommand{\psl}{PSL_2 (\F_p)}
\newcommand{\psls}{PSL_2 (\F_7)}
\newcommand{\F}{\mathcal{F}}
\newcommand{\injects}{\hookrightarrow}
\newcommand{\vect}[1]{\stackrel{\rightharpoonup}{\mathbf #1}}
\newcommand{\Css}{\mathcal{C}_{ss}}

\newcommand{\G}{\mathcal{G}}
\newcommand{\E}{\mathcal{E}}
\newcommand{\V}{\mathcal{V}}
\newcommand{\mL}{\mathcal{L}}
\newcommand{\D}{\mathcal{D}}
\newcommand{\I}{\mathcal{I}}
\newcommand{\mH}{\mathcal{H}}
\def\endrem{}
\def\colon{{:}\;}
\newcommand{\Rep}{\mathrm{Rep}}
\newcommand{\Hom}{\mathrm{Hom}}
\newcommand{\Lie}{\mathrm{Lie}}
\newcommand{\K}{K_{\mathrm{def}}}
\newcommand{\mO}{\mathcal{O}}
\newcommand{\Map}{\mathrm{Map}}
\newcommand{\flatc}{\mathcal{A}_{\mathrm{flat}}}
\newcommand{\A}{\mathcal{A}}
\newcommand{\hofib}{\mathrm{hofib}}
\newcommand{\Cmin}{\mathcal{C}_{\textrm{min}}}

\newcommand{\cA}{\mathcal{A}}
\newcommand{\cG}{\mathcal{G}}
\newcommand{\fg}{\mathfrak{g}}
\newcommand{\ad}{\mathrm{ad}}
\newcommand{\Ad}{\mathrm{Ad}}
\newcommand{\Aut}{\mathrm{Aut}}
\newcommand{\wt}[1]{\widetilde{#1}}
\newcommand{\expl}{\mathrm{exp}}

\title{The Yang--Mills stratification for surfaces revisited}

\author{Daniel A. Ramras}

\address{Department of Mathematical Sciences,
New Mexico State University, Las Cruces, New Mexico 88003-8001, USA}
\email{ramras@nmsu.edu}

\thanks{This work was partially supported by
    NSF grants DMS-0353640 (RTG), DMS-0804553, and DMS-0968766}

\begin{abstract}
We revisit Atiyah and Bott's study of Morse theory for the Yang--Mills functional over a Riemann surface, and establish new formulas for the minimum codimension of a (non-semi-stable) stratum.  These results yield the exact connectivity of the natural map 
$$(\Cmin(E))_{h\G(E)}\to \Map^E (M, BU(n))$$
from the homotopy orbits of the space 
of central Yang--Mills connections to the classifying space of the gauge group $\G(E)$.  All of these results carry over to non-orientable surfaces via Ho and Liu's non-orientable Yang--Mills theory.   

A somewhat less detailed version of this paper (titled ``On the Yang--Mills stratification for surfaces") will appear in the Proceedings of the AMS.
\end{abstract}

\maketitle{}

\section{Introduction}

Let $M^g$ be a Riemann surface of genus $g>0$, and consider a vector bundle $E$ over $M$.  When $E$ is trivial, the space $\flatc(E)$ of flat connections on $E$ forms the minimum critical set for the Yang--Mills functional $L:\cA(E) \to \bbR$, where $\cA(E)$ is the affine space of all connections, and for $A\in \cA(E)$, 
$$L(A) = \int_M ||F(A)||^2 d\mathrm{vol}_M.$$
Here $F(A)$ is the curvature form of $A$ and the volume of $M$ is normalized to 1.
In their seminal paper on Yang--Mills theory~\cite{A-B}, Atiyah and Bott showed that by treating $L$ as a gauge-equivariant Morse function, one can learn a great deal about the topology of the critical set $\flatc(E)$ and its stable manifold $\C_{\mathrm{ss}} (E)$, the space of semi-stable holomorphic structures on $E$.  In particular, Atiyah and Bott provided a framework for computing the gauge-equivariant cohomology of these spaces.  
When $E$ is non-trivial, the minimum critical  set $\Cmin (E)$ of the Yang--Mills functional consists of \emph{central} Yang--Mills connections and similar methods may be used to study the topology of this space.  Again the space  $\C_{\mathrm{ss}} (E)$ of semi-stable structures serves as the stable manifold of this critical set.  

This theory has been fleshed out and made rigorous by the work of Daskalopoulos~\cite{Dask} and R{\aa}de~\cite{Rade} (building upon Uhlenbeck's famous compactness theorem).  In this paper, we prove (Theorem~\ref{conn-calc})
that if  $E\to M$ has rank $n$ and Chern number $k$, then the connectivity of $\Cmin(E)$ is precisely $2\gcd (n,k) - 2$ if the genus $g$ of $M$ is 1, and precisely $2\left(\min ([k]_n, [-k]_n) + (g-1)(n-1)\right) - 2$ if $g>1$ (here $[r]_n$ denotes the unique integer between $1$ and $n$ congruent to $r$ modulo $n$).  
We provide a similarly explicit formula (Theorem~\ref{conn-calc-non-orient}) in the case of a non-orientable surface $M$, using Ho and Liu's non-orientable Yang--Mills theory~\cite{Ho-Liu-non-orient}.  The same formulas also give the connectivity
of the natural map
$$E\G(E) \cross_{\G(E)} \Cmin (E) = \Cmin (E)_{h\G(E)}\to \Map^{E} (M, BU(n))$$
from the homotopy orbits of $\Cmin(E)$ under the gauge group $\G(E)$ to the classifying space of 
$\G(E)$.  

These results rely on an interesting combinatorial analysis of the Yang--Mills stratification.   In the orientable case, weaker bounds on this quantity have been used in the literature before (see Daskalopoulos \cite[Section 7]{Dask}, Daskalopoulos and Uhlenbeck~\cite{Dask-Uhl}, Ramras~\cite{Ramras-surface}, and Cohen--Galatius--Kitchloo~\cite{CGK}), but our precise formulas are new.  Our connectivity results also provide isomorphisms 
$$H^*_{\G(E)} (\Cmin (E)) \isom H^*(B\G(E); \bbZ)$$
and
$$H^*_{\G_0(E)} (\Cmin(E)) = H^*(\Cmin(E)/\G_0(E))\isom H^*(B\G_0(E);\bbZ)$$ 
in low dimensions.
In the orientable case, $H^*(B\G(E); \bbZ)$ and $H^*(B\G_0(E);\bbZ)$ were computed in Atiyah--Bott~\cite[Section 2]{A-B}.  For non-orientable surfaces, rational cohomology may be computed similarly (see Ho--Liu~\cite[Section 2]{Ho-Liu-or-non}).

Atiyah and Bott's approach to calculating the gauge-equivariant cohomology of $\Cmin(E)$ may be seen as a close analogue of finite dimensional Morse theory. 
In the finite dimensional setting, one describes a manifold as a cell complex built inductively according to the critical-point structure of the Morse function in question.  Each critical point corresponds to the addition of a new cell whose dimension equals the index of that critical point, i.e. the codimension of the stable manifold.  In the infinite dimensional setting of Yang--Mills theory, one tries to mimic this picture by building up the space of connections one Yang--Mills stratum $\C_\mu$ at a time.
At each stage, rather than attaching a finite dimensional cell, we add a new finite codimension submanifold.  The effect in gauge-equivariant (co)homology can be analyzed by considering the long
exact sequence associated to the pair $(\bigcup_{i=1}^{n+1} \C_{\mu_i} , \bigcup_{i=1}^{n} \C_{\mu_i})$.  One then hopes to establish a Thom isomorphism for the relative terms:
$$H^*_{\G(E)} \left(\bigcup_{i=1}^{n+1} \C_{\mu_i},  \bigcup_{i=1}^{n} \C_{\mu_i}; \bbZ\right) \isom H^{*-\mathrm{codim} (\C_{\mu_m})}_{\G(E)} \left(\C_{\mu_m}; \bbZ \right).$$   
Atiyah and Bott used this method to calculate the gauge-equivariant integral cohomology of the space of semi-stable holomorphic structures. 

However, this approach relies on several technical points:
first, the strata must be locally closed submanifolds, of finite codimension in the space of all connections.
This issue was resolved by Daskalopoulos~\cite{Dask}.
Next, one must establish the necessary Thom isomorphisms.  In Ramras \cite{Ramras-tubular}, it is shown that
the strata $\C_{\mu}$ have gauge-invariant tubular neighborhoods (neighborhoods equivariantly homeomorphic the \emph{orientable} normal bundles $\nu (\C_{\mu})$), and then the desired Thom isomorphisms follow by excising the complements of these tubular neighborhoods and applying the ordinary Thom isomorphism theorem (see Corollary~\ref{Thom}).  An alternate approach would be to use the fact that after modding out the (based, complex) gauge group, one obtains a smooth algebraic variety, and results of Shatz~\cite[Section 4]{Shatz} show that the bundles of type $\mu$ form a smooth subvariety.  This algebraic approach involves various technicalities, which we will not attempt to resolve here.  Our infinite dimensional approach also yields tubular neighborhoods (and Thom isomorphisms) in the setting of Ho and Liu's non-orientable Yang--Mills theory.

 Orientability of the normal bundles is immediate over Riemann surfaces, because these are naturally \emph{complex} vector bundles.  Over a non-orientable surface, these are only \emph{real} vector bundles, and before applying the (integral) Thom isomorphism theorem it is necessary to know that these bundles are orientable.  A solution to this problem is given in Ho--Liu--Ramras~\cite{Ho-Liu-Ramras} 
 
 For the connectivity calculations in this paper, the mod--2 Thom isomorphism in non-equivariant cohomology, together with the universal coefficient theorem suffice.  Hence these results do not depend on the existence of gauge-invariant tubular neighborhoods.  There are still some subtleties regarding the construction of tubular neighborhoods for non-closed submanifolds, and we discuss these issues in Section~\ref{inv}.

This paper is organized as follows.  In Section~\ref{hom-stab}, we explain Atiyah and Bott's inductive method for calculating gauge-equivariant (co)homology using the Yang--Mills stratification, and establish the necessary combinatorial properties of this stratification (filling in some gaps in the literature).
In Section~\ref{inv}, we discuss tubular neighborhoods for locally closed submanifolds, and in Section~\ref{conn} we combine these results with a combinatorial analysis of the Yang--Mills stratification to obtain the connectivity calculations described above.

A somewhat less detailed version of this paper (titled ``On the Yang--Mills stratification for surfaces") will appear in the Proceedings of the AMS.

\vspace{.1in}

\noindent {\bf Acknowledgements:} I thank G. Carlsson, C. Groft, G. Helleloid, and C.-C. Liu for helpful conversations.  Additionally, I thank N.-K. Ho for pointing out several misstatements in an earlier draft.

\section{The Harder--Narasimhan stratification}$\label{hom-stab}$

In this section we recall and analyze the Harder--Narasimhan stratification on the space of holomorphic structures on a smooth, complex vector bundle over a Riemann surface $M=M^g$, as in~\cite[Section 7]{A-B} (we suppress the genus $g$ when possible).  This stratification agrees with the Morse stratification for the Yang--Mills functional, in the sense that the Yang--Mills flow defines deformation retractions from each Harder--Narasimhan stratum to its subset of Yang--Mills critical points~\cite{Dask, Rade}.  

Let $\C(E) = \C(n,k)$ denote the space of holomorphic structures on a rank $n$ Hermitian bundle $E$ with Chern number $k$. 
As shown in Atiyah--Bott~\cite[Sections 5, 7]{A-B}, this is an affine space, isomorphic to the affine space 
$\A(E)$ of Hermitian connections on $E$.  As such, we may equip this space with a Sobolev norm and complete it to a Hilbert space.  Throughout this paper, $\C(n,k)$ will denote such a Hilbert space completion (we will not need to specify the Sobolev regularity).  Recall that the unitary gauge group $\G(E) = \G(n,k)$ of unitary automorphisms of $E$, and the larger \emph{complex} gauge group $\G^\bbC (E) = \G^\bbC (n,k)$ of all complex automorphisms of $E$, act on the space $\C(n,k)$.  We will always implicitly consider the Hilbert Lie-group completions of these groups, as in~\cite[Section 14]{A-B}.

To define the Harder--Narasimhan stratification of $\C(n,k)$, we must first recall the Harder--Narasimhan filtration on a holomorphic bundle.
Given a holomorphic structure $\E$ on a bundle $E \to M$  of rank $n$ and Chern number $k$, there is a unique filtration (the Harder--Narasimhan filtration)
$$0 = \E_0 \subset \E_1 \subset \cdots \E_r = \E$$
of $\E$ by holomorphic subbundles with the property that each quotient $\D_i = \E_i/\E_{i-1}$ is semi-stable ($i = 1,\ldots ,r$) and $\mu(\D_1) > \mu(\D_2) > \cdots > \mu(\D_r)$, where the ``slope" $\mu(D_i)$
is defined by $\mu(D_i) = \frac{\mathrm{deg} (D_i)}{\mathrm{rank} (D_i)}$.
(Recall that a bundle $F$ is semi-stable if for all holomorphic subbundles $F'<F$, $\mu(F')\leqs \mu(F)$.)

Letting $n_i =$ rank$(D_i)$ and $k_i =$ deg$(D_i)$, we call the sequence 
$$\mu = ((n_1, k_1), \ldots, (n_r, k_r))$$
the \emph{type} of $\E$.  Let $\C_{\mu} = \C_{\mu}(n,k) \subset \C(n,k)$ denote the subspace of all holomorphic structures \emph{complex} gauge-equivalent to a smooth structure of type $\mu$ (by Atiyah--Bott~\cite[Section 14]{A-B}, every orbit of the Hilbert Lie group $\G^\bbC (n,k)$  on the Hilbert space $\C(n,k)$ contains a unique isomorphism class of holomorphic structures).  Note that the semi-stable stratum corresponds to $\mu = ((n,k))$, and that since degrees add in exact sequences we have
$\sum_i k_i = k$.
With this notation, we now have the following result from~\cite[Section 7]{A-B} (see also~\cite[Theorem B]{Dask}).

\begin{theorem}$\label{strata}$
Let $\mu = ((n_1, k_1), \ldots, (n_r, k_r))\in \C(n,k)$.  Then the stratum $\C_{\mu}$ is a locally closed submanifold of $\C(n,k)$
with complex codimension given by
$$c(\mu) = \sum_{i>j} n_i k_j - n_j k_i + n_i n_j (g-1).$$
\end{theorem}

Following Atiyah and Bott, we proceed to describe $\C(n,k)$ as a colimit over unions of strata.  This facilitates the inductive calculation of equivariant homology from~\cite{A-B}.  
The following definition will be useful.

\begin{definition}$\label{admissible}$ A sequence $((n_1, k_1), \ldots, (n_r, k_r))$ is \emph{admissible} of total rank $n$ and total degree $k$ if 
$n_i > 0$ for each $i$, $\sum n_i = n$, $\sum k_i = k$, and $\frac{k_1}{n_1}>\cdots>\frac{k_r}{n_r}$.    
We denote the set of all admissible sequences of total rank $n$ and total degree $k$ by $\mathcal{I}(n, k)$.  
The set $\I(n,k)$ has a partial ordering defined as follows: given an admissible sequence $\mu = ((n_1, k_1), \ldots, (n_r, k_r))$, let $\hat{\mu} = (\hat{\mu}_1,\hat{\mu}_2,\ldots, \hat{\mu}_n)$ where the first $n_1$ terms equal $k_1/n_1$, and next $n_2$ equal $k_2/n_2$ and so on.  Then we say $\lambda\geqs \mu$ if 
$$\sum_{j\leqs i} \hat{\lambda}_j \geqs \sum_{j\leqs i} \hat{\mu}_j$$
for $i = 1,\ldots, n$.
\end{definition}

Following~\cite{A-B}, we introduce another way of thinking about the ordering on these strata (due to Shatz~\cite{Shatz}).  Given an admissible sequence $\mu$, we construct a convex path $P(\mu)$ in $\bbR^2$ starting at $(0,0)$ and ending at $(n,k)$ by connecting the points 
$(\sum_{j=1}^i n_j, \sum_{j=1}^i k_j)$ with straight lines ($i = 1,2,\cdots n$).  Convexity corresponds precisely to the condition that the slopes decrease, i.e. that 
$$\frac{k_1}{n_1}>\frac{k_2}{n_2}>\cdots>\frac{k_r}{n_r}.$$
Now, for any $\lambda, \mu\in \I(n,k)$, we have
$\lambda \geqs \mu$ if and only if $P(\lambda)$ lies above $P(\mu)$.  Note that we may recover the sequence $\mu$ from $P = P(\mu)$ by reading off the coordinates of the points where $P$ changes slope, and any convex path from $(0,0)$ to $(n,k)$ which changes slope only at points with integer coordinates yields an admissible sequence.

\begin{remark}$\label{genus0}$
When $g=0$, Grothendieck's theorem states that every holomorphic bundle is a sum of line bundles.  Hence in genus zero, the stratum corresponding to $\mu\in \I(n,k)$ may empty.  For this reason, we assume $g>0$ throughout this paper (and in the non-orientable case we do not consider $\bbR P^2$).
\end{remark}

The necessary fact regarding the Harder--Narasimhan stratification is the following result, essentially due to Atiyah and Bott.  Here we will fill in some details of the proof absent from their paper~\cite{A-B}, and which do not appear to have been clarified in the literature.

\begin{proposition}$\label{total-order}$ The partial ordering $\leqs$ on $\I(n,k)$ can be refined to a linear ordering $\mu_1\prec \mu_2 \prec \cdots$  such that
for any $j$, $\C_j = \C_j (n,k) = \bigcup_{i=1}^j \C_{\mu_i}$ is open in $\C(n,k)$.  
\end{proposition}

Let $E\to M$ be a Hermitian bundle over a non-orientable surface, and let $\widetilde{E}\to \widetilde{M}$ denote the pullback of $E$ to the orientable double cover of $M$.  Connections on $E$ pull back to connections on $\widetilde{E}$, yielding an embedding $i:\A(E)\injects \A(\widetilde{E})$,
and following Ho and Liu~\cite{Ho-Liu-non-orient} we define the Yang--Mills strata of $\A(E)$ to be
the intersections of $\A(E)$ with the Harder--Narasimhan strata of $\A(\widetilde{E})\isom \C(\widetilde{E})$.  Now Proposition~\ref{total-order} implies:

\begin{corollary}$\label{total-order-non-orient}$
For any Hermitian bundle $E$ on a non-orientable surface, 
the linear ordering $\prec$ on $\A(\widetilde{E})$ induces a linear ordering on the Yang--Mills strata of $\A(E)$ such that the union of any initial segment $\{S| S\prec S_0\}$ is open in $\A(E)$.
\end{corollary}

\begin{remark}$\label{empty-strata}$
Say $\Sigma$ is a non-orientable surface with orientable double cover $M^g$, and $E\to \Sigma$ is a Hermitian bundle.  Then the intersections of the Harder-Narasimhan strata for the pullback $\wt{E} \to M^g$ with the space of connections on $E$ are sometimes empty.  See~\cite[Section 7.1]{Ho-Liu-non-orient} for a precise determination of the non-empty strata.
\end{remark}

%\vspace{.15in}
The proof of Proposition~\ref{total-order} will require several lemmas, all implicit in~\cite{A-B}.  First we note some simple but important corollaries.  One would like to compute (co)homology inductively, by analyzing the spectral sequence (i.e. the collection of long exact sequences) associated to the filtration of $\C(n,k)$ by the strata.  At each stage, one wants to analyze the relative term
$H^*(C_{m}, C_{m-1})$.  The key result is the following Thom isomorphism.

\begin{corollary}$\label{Thom}$  Let $M^g$ be a Riemann surface of genus $g>0$, and let $E$ be a complex vector bundle over $M$.  Let $\C_{\mu_1}\prec \C_{\mu_2} \prec \cdots$ be a linear order on the Harder--Narasimhan strata of $\C(E)$ as in Proposition~\ref{total-order}.  Then there are Thom isomorphisms 
$$H_* (\C_m, \C_{m-1}; \bbZ) \isom H_{*-\mathrm{codim} (\C_{\mu_m})} (\C_{\mu_m}; \bbZ)$$
and similarly for integral cohomology.  Moreover, analogous isomorphisms hold for integral gauge-equivariant homology and cohomology.  Here $H_*$ and $H^*$ are interpreted as zero when $*$ is negative.

The same results hold in the space of connections on a complex bundle over any non-orientable surface $\Sigma$, so long as the genus $\wt{g}$ of the orientable double cover $\wt{\Sigma}$ is at least 2.  With $\bbZ/2\bbZ$--coefficients, these results hold even when $\wt{g} = 1$.
\end{corollary}
\noindent {\bf Proof.}  In the orientable case, this corollary is a simple consequence of the results proven in this paper.  In Ramras \cite{Ramras-tubular}, we construct a $\G(E)$--invariant tubular neighborhood $\nu_m\subset  \C_m$ of the not-necessarily-closed submanifold $\C_{\mu_m}$.  Excising the complement of $\nu_m$ in $\C_m$ and applying the Thom isomorphism theorem to the (complex, hence orientable) normal bundle of $\C_{\mu_m}$ gives the desired isomorphism in the non-equivariant case.  In the equivariant case, one simply observes that for any $G$--equivariant complex vector bundle $V\to X$, the homotopy orbit bundle $V_{hG}\to X_{hG}$ is still a complex vector bundle.
Hence we may excise the complement of $(\nu_m)_{h\G(E)}$ in $(\C_m)_{h\G(E)}$ and apply the ordinary Thom isomorphism to the bundle $(\nu_m)_{h\G(E)}\to (\C_{\mu_m})_{h\G(E)}$.

The same argument works in the non-orientable case, since it is proven in Ho--Liu--Ramras~\cite{Ho-Liu-Ramras} that the normal bundles to the Yang--Mills strata (and the corresponding homotopy orbit bundles) are orientable (real) vector bundles, so long as $\wt{g}\geqs 2$.
$\hfill \Box$  \vspace{.35in}

We now explain how to compute (gauge-equivariant) (co)homology of $\C(E)$ inductively via the linear ordering on the set of Harder--Narasimhan strata.
For the cohomological case, we need a simple finiteness property of this stratification.
As observed by Atiyah and Bott~\cite[p. 569]{A-B}, this lemma follows quickly from Theorem~\ref{strata}; for completeness we provide a proof.

\begin{lemma}[Atiyah-Bott]$\label{small-codim}$
For any $n, D\in \bbN$ and any $k\in \bbZ$, there are finitely many 
admissible sequences $\mu\in \I(n,k)$ with $c(\mu)<D$.
\end{lemma}

\noindent {\bf Proof.}  Let $\mu = ((n_1, k_1), \ldots, (n_r, k_r))$ be an admissible sequence with
$c(\mu)< D$.  Since $\sum n_i = n$, there are finitely many possibilities for the positive integers $n_i$.  
By convexity, we have $k_1/n_1 > k/n$, and hence $k_1 > \frac{kn_1}{n}$.  When $k\geqs 0$, this means that $k_1>k/n$; when $k<0$ it means that $k_1>k$.
We will check that if $k\geqs 0$, then $k_i > -D$ for each $i>1$, and if $k<0$, then $k_i> k-D$ for $i>1$; since $\sum k_i = k$ this means there are finitely many possibilities for the integers $k_i$.

Since each term
in the sum defining $c(\mu)$ is positive (Theorem~\ref{strata}) we know that 
$k_1 n_i - k_i n_1<D$ for each $i$, and rearranging gives
$k_i > \frac{k_1 n_i - D}{n_1}$.
 We now use our bounds on $k_1$. 
When $k\geqs 0$, we have $\frac{ k_1 n_i- D}{n_1}>\frac{(k/n) n_i  - D}{n_1} \geqs \frac{-D}{n_1}\geqs -D$ as desired.  When $k<0$, we have
$\frac{ k_1 n_i- D}{n_1} > \frac{(kn_1/n) n_i  - D}{n_1} > k-D$.
$\hfill \Box$  \vspace{.35in}

\begin{corollary}$\label{gauge-hom}$ For any Hermitian bundle $E$ over a Riemann surface, there are isomorphisms
$$H_* (B\G(E)) \isom  H^{h\G(E)}_* (\C(E); \Z) \isom 
     \colim_{j\to \infty} H^{h\G(E)}_* \left(C_j (E); \Z \right)$$
and
$$H_* (B\G(E)) \isom H_{h\G(E)}^* (\C(E); \Z) \isom 
     \lim_{\stackrel{\leftarrow}{j}} H_{h\G(E)}^* \left(C_j (E); \Z \right)$$
in gauge-equivariant integral (co)homology.
(The corresponding statements for ordinary (co)homology with arbitrary coefficients also hold.)

If $E$ is a Hermitian bundle over a non-orientable surface $\Sigma$, then the same statements hold for the flitration of $\A(E)$ induced by the stratification in Corollary~\ref{total-order-non-orient} 
(although for equivariant cohomology we must assume that the genus of the orientable double cover $\wt{\Sigma}$ is greater than 1).
\end{corollary}
\noindent {\bf Proof.}  The left-hand isomorphisms follow from contractibility of the affine space $\C(E)$ (or, in the non-orientable case, $\A(E)$).  Proposition~\ref{total-order} and Corollary~\ref{total-order-non-orient}
immediately yield the right-hand isomorphisms at the level of (equivariant) singular chains and cochains.  The homological results follow from the fact that homology commutes with directed limits.  In cohomology, Lemma~\ref{small-codim} and Corollary~\ref{Thom} imply that for each $p$ the inverse system $\{H^p_{\G(n,k)} (\C_j(n,k); \Z)\}_j$ is eventually constant, so $\lim^1$ vanishes and the result follows from~\cite{Milnor-lim1}.
$\hfill \Box$  \vspace{.15in}

The proof of Proposition~\ref{total-order} will require
one further finiteness property of the partial ordering on $\I(n,k)$, also noted by Atiyah and Bott~\cite[p. 567]{A-B}.  

\begin{lemma}[Atiyah--Bott]$\label{finiteness}$
If $I\subset \mathcal{I}(n, k)$ is a finite collection of admissible sequences, then there are finitely many minimal elements in the complement $I^c = \mathcal{I}(n, k) - I$.  
\end{lemma}
\noindent {\bf Proof.}  We will phrase the argument in terms of convex paths.  Let $I$ be a finite collection of convex paths from $(0,0)$ to $(n,k)$.  If $P$ is a minimal path in the complement of $I$, then every path beneath $P$ lies in $I$, so either $P$ is the minimum path, i.e. the line from $(0,0)$ to $(n,k)$, or $P$ is a minimal cover of a path $Q\in I$, meaning that $Q<P$ and there is no path $P'$ with $Q<P'<P$.  So to prove the first statement of the lemma, it will suffice to show that each path $Q$ has only finitely many minimal covers.  In the course of proving this fact, we will also prove the second statement of the lemma.

Fix a sequence  $\mu = ((n_1, k_1), \ldots, (n_r, k_r))\in I(n,k)$ and let $P = P(\mu)$ be the associated path.  
Define
$$s_1 (P) = \max \{k_1/n_1, 0 \}; \,\,\,\, s_r (P) = \min \{k_r/n_r, 0\}.$$
Consider another path $Q = P(\nu)$, where $\nu \in I(n,k)$ and $\nu\neq ((n,k))$.  
Let $h(Q) = (h_1 (Q), h_2 (Q))$ denote the right endpoint of the rightmost line segment in $Q$ with slope at least $\frac{k}{n}$.  We claim that if $h_2 (Q) \geqs n(s_1(P) - s_r (P)) + \max \{k, 0\} + 1$, and $Q'$ has vertices 
$(0,0)$, $(h_1 (Q), h_2 (Q) - 1)$ and $(n, k)$, then
\begin{equation}\label{Q'}
P\leqs Q'<Q \textrm{\,\,\,\, and if \,\,} r > 2 \mathrm{\,\, then\,\,} P<Q<Q'
\end{equation}
Assuming (\ref{Q'}), we now complete the proof.  If $Q$ is a minimal cover of $P$ 
then either $r>2$ and $h_2 (Q) \leqs n(s_1(P) - s_r (P)) + \max \{k, 0\}$, or $r=2$ and $P=Q'$.  In the former case, $Q$ lies below the line of slope $\frac{k}{n}$ passing though the point $h(Q) = (h_1 (Q), h_2 (Q))$.  Since $h_1 (Q) \leqs n$, this restricts $Q$ to a finite region (depending only on $P$).
The latter conditions can hold for only finitely many paths $Q$, since if $P = Q'$ then $Q$ passes through $(n, k_1 + 1)$.

To prove (\ref{Q'}), first note that 
$Q'< ((0,0), (h_1 (Q), h_2 (Q)), (n,k)) \leqs Q$,  
so we need only check that $P<Q'$ when $r>2$.  If not, then at some time $x=x_0$ the path $P$ lies above the path $Q'$ (since $r>2$, $P\neq Q'$).  If $x_0 \leqs h_1 (Q)$,
then since $Q'$ is just a straight line for $x\leqs h_1 (Q)$, the initial slope of $P$ must be more than the initial slope of $Q'$.  Our assumption on $h_2 (Q)$ now gives
\begin{equation}\label{s_1}
s_1 (P) \geqs \frac{k_1}{n_1} > \frac{h_2 (Q) - 1}{h_1 (Q)} \geqs \frac{(n(s_1 (P) - s_r (P)) + \max \{k, 0\}  +1) - 1}{n}.
\end{equation}
Since $s_r(P) \leqs 0$, (\ref{s_1}) yields
\begin{equation*}
\begin{split}
s_1 (P) &>  \frac{(n(s_1 (P) - s_r (P)) + \max \{k, 0\}  +1) - 1}{n} \geqs \frac{n(s_1 (P) - s_r (P))}{n}\\
& = s_1(P) - s_r(P) \geqs s_1(P),
\end{split}
\end{equation*}
a contradiction.
Similarly, if $x_0 > h_1(Q)$ then the final slope of $P$ is less than the final slope of $Q'$.  Moreover, $s_1(P)\geqs 0$ and $s_r(P)\leqs 0$, so we have
\begin{equation*}
\begin{split} 
s_r (P) & \leqs \frac{k_r}{n_r} < \frac{k-(h_2(Q) - 1)}{n - h_1 (Q)} \leqs \frac{k-(n(s_1(P) - s_r (P)) + \max\{k,0\})}{n - h_1 (Q)}\\
&\leqs \frac{-ns_1(P) + n s_r (P)}{n - h_1 (Q)}  \leqs \frac{n s_r(P)}{n - h_1 (Q)} = \frac{n}{n-h_1(Q)} s_r(P) 
\leqs s_r (P),
\end{split}
\end{equation*}
a contradiciton as before.
$\hfill \Box$  \vspace{.15in}

We note that the proof of Lemma~\ref{finiteness} also shows that for any $\mu\in I(n,k)$, all but finitely many $\lambda\in I(n,k)$ satisfy $\lambda\geqs \mu$.
The final ingredient in the proof of Proposition~\ref{total-order} is the following result regarding the closures on the Harder--Narasimhan strata.

\begin{proposition}$\label{closures}$
Let $S\subset \I(n,k)$ be a collection of admissible sequences that is upwardly closed, in the sense that if $\mu > \mu'$ and $\mu'\in S$, then $\mu\in S$ as well.
Then 
the set $\bigcup_{\mu \in S} \C_{\mu}$ is closed.
\end{proposition}

 Atiyah and Bott~\cite[(7.8)]{A-B}, as well as Daskalopoulos~\cite[Proposition 2.12]{Dask}, state only the (strictly weaker) fact
\begin{equation}\label{AB-closure-1}
\overline{\C_{\mu}} \subset \bigcup_{\mu'\geqs \mu} \C_{\mu'},
\end{equation}
where $\overline{\C_\mu}$ denotes the closure of this stratum (this result originated in the algebro-geometric work of Shatz~\cite{Shatz}).  To prove the stronger statement in Proposition~\ref{closures}, we will apply another result of Atiyah and Bott~\cite[Section 8]{A-B}.

\begin{proposition}[Atiyah-Bott]\label{inf}  
Consider an admissible sequence 
$$\mu = ((n_1, k_1), \ldots, (n_r, k_r))\in I(n,k).$$  
Then for any $A\in \C_{\mu}$, we have
$$l(\mu) := \inf_{g\in \G^{\bbC}(n,k)} L(g\cdot A) = \sum_{i=1}^{r} \frac{k_i^2}{n_i}$$
where $L$ denotes the Yang--Mills functional and the infimum is taken over the \emph{complex} gauge group.
\end{proposition}

\noindent {\bf Proof of Proposition~\ref{closures}.} By (\ref{AB-closure-1}), we have $\bigcup_{\mu\in S} \C_{\mu} = \bigcup_{\mu \in S} \overline{\C_{\mu}}$.  Since the union of a locally finite collection of closed sets is closed, it will suffice to show that the closures of the strata $\C_\mu$ form a locally finite cover of $\C(n,k)$.  We will check that for each $N\in \bbZ$, only finitely many closures $\overline{\C_{\mu}}$ contain elements $A$ with $L(A) <N$.  

For any $M\in \bbR$ there are finitely many $\mu\in \mathcal{I}(n, k)$ with $l(\mu)\leqs M$, because $l(\mu)\leqs M$ implies that the path $P(\mu)$ lies entirely under the line
$y=\sqrt{M}x$.  It now suffices to check that if $L(A)<N$ for some $A\in \overline{\C_{\mu}}$, then
$l(\mu)<N$.  By continuity of $L$, there exists $A'\in \C_\mu$ with $L(A')<N$, and Proposition~\ref{inf} implies that $l(\mu) \leqs L(A')$.
$\hfill \Box$

\begin{remark}$\label{l(mu)}$
Although we will not need this fact, we point out that the number $l(\mu)$ appearing in Proposition~\ref{inf} is actually the (unique) critical value of the Yang--Mills functional on the stratum $\C_\mu$.  This follows from convergence of the Yang--Mills flow (R{\aa}de~\cite{Rade}) and the fact that the Morse strata agree with the Harder-Narasimhan strata (Daskalopoulos \cite{Dask}), together with discreteness of the critical values of $L$.  As mentioned in R{\aa}de \cite[Section 2]{Rade}, this follows from Uhlenbeck Compactness and  \cite[Proposition 7.2]{Rade}.
\end{remark}

%\vspace{.15in}

\noindent {\bf Proof of Proposition~\ref{total-order}.}  
We construct a linear ordering $\prec$ on $\I(n,k)$ by setting $T_0 = \{ ((n,k))\}$, and inductively defining
$$T_l = T_{l-1} \cup \{\mu\in \I(n,k) \,\,|\,\, \mu \textrm{\,\, is minimal in \,\,} \I(n,k)\setminus T_{l-1} \},$$
where we choose any linear ordering $\prec$ on $T_l$ which extends the existing ordering $\prec$ on $T_{l-1}$ and satisfies $\mu \prec \eta$ if $\mu\in T_{l-1}$ and $\eta \in T_l \setminus T_{l-1}$.
The set $T = \bigcup_l T_l$ is linearly ordered by $\prec$, and if $\mu\leqs \eta$ then $\mu\prec \eta$.  

We must check that $T = \I(n,k)$.  If $\I(n,k)\setminus T$ is non-empty, then we may choose a $\leqs$--minimal element $\mu$ from this poset.
There are finitely many $\eta$ with $\eta< \mu$, so we may choose a $\leqs$--maximal element $\eta_0$ from the finite set 
$T_\mu = T\cap \{\eta | \eta < \mu\}$ (note that $T_\mu$ necessarily contains the minimum sequence $((n,k))$, so $T_\mu\neq \emptyset$).  Then $\eta_0\in T_N$ for some $N$.  If $\mu$ were a minimal cover of $\eta_0$ in the poset $(\I(n,k), \leqs)$, then by definition we would have $\mu\in T_{N+1}$, a contradiction.  So we may choose a minimal cover $\eta_1$ of $\eta_0$ with $\eta_0 < \eta_1 < \mu$.  Then $\eta_1\in T_{N+1}$, so $\eta_1 \in T_\mu$, contradicting maximality of $\eta_0$.  Hence $\I(n,k) \setminus T$ must be empty.
$\hfill \Box$

%%%%%%%%%%%%%%%%%%%%%%%%%%%%%%%%%%%%%%%%%%%%%%%%%%%
%%%%%%%%%%%%%%%%%%%%%%%%%%%%%%%%%%%%%%%%%%%%%%%%%%%

\section{Tubular Neighborhoods}$\label{inv}$

In the next section, we will need to apply the Thom isomorphism theorem to calculate $H^*(\C_i, \C_{i-1}; \bbZ/2\bbZ)$.  This depends on the existence of tubular neighborhoods for the Yang--Mills strata $\C_{\mu_i}$ inside the open sets $\C_i =  \bigcup_{i=1}^j \C_{\mu_i}$.  Although the construction is by now more or less standard, some subtleties arise due to the fact that the sets $\C_{\mu_i}$ are not closed.  Hence we outline the argument.  Our basic reference for Banach manifolds is Lang~\cite{Lang-dg}.

\begin{lemma}$\label{tubular}$
Let $Y$ be a smooth, metrizable Banach manifold, and let $X\subset Y$ be a locally closed submanifold of finite codimension.  Then there is an open neighborhood $\tau(X)$ of $X$ inside $Y$ which is diffeomorphic to the normal bundle $N(X)$.
\end{lemma}
\noindent {\bf Proof.} We follow Lang~\cite[Section IV.5]{Lang-dg} and Bredon~\cite[VI.2]{Bredon-transf}.  As shown in~\cite[Chapter III]{Lang-dg}, there is a smooth direct sum decomposition $T(Y)|_X = T(X) \oplus N(X)$.  Using sprays, Lang constructs an exponential map $\expl: \D \to Y$, where $\D\subset T(Y)$ is an open neighborhood of the zero section, and shows that $\expl$ restricts to a local diffeomorphism $\D\cap N(X) \to Y$ (this means each $x\in X$ has an open neighborhood $U_x \subset N(X)$ on which $\expl$ is a diffeomorphism onto an open set in $Y$).  

We claim that there exists an open set $W\subset \D\cap N(X)$ with the property that $\{w\in W \, | \, \expl(w) \in X\} = W\cap X$.
Since $U_x \cap X$ is open in $X$, we have $U_x \cap X = W_x \cap X$ for some open set  $W_x \subset Y$ .  Now  
$W = \bigcup_x \left(U_x\cap  \expl^{-1} (W_x) \right)$ is the desired open set in $N(X)$.
Bredon~\cite[Chapter VI, Lemma 2.3]{Bredon-transf} now shows that there exists a smaller neighborhood $W'\subset W$ on which $\expl$ is injective (Bredon  assumes  both $W$ and $Y$ are metric spaces, but only uses the fact that $Y$ is  metric).

Now $\expl: W' \to \expl(W')$ is a diffeomorphism onto an open neighborhood of $X$ inside $Y$.  Finally, Lang~\cite[Section VII.4]{Lang-dg} shows that the finite-dimensional vector bundle $N(X)$ can be ``compressed" into the neighborhood $W'$.
$\hfill \Box$  

\begin{remark} Lang~\cite[Section IV.5]{Lang-dg} assumes that $X$ is closed in $Y$.  This is used in the construction of the neighborhood $W'$ above.  One must replace a family $\{U_\alpha\}_\alpha$ of open sets in $Y$ which covers $X$ by a subordinate family $\{V_\beta\}_\beta$ which still covers $X$ and for which $\{\overline{V_\beta}\}_\beta$ is locally finite.  This can always be done if $X$ is closed and $Y$ admits partitions of unity (Lang's hypotheses) but also works if $Y$ is metrizable, because then the union of the $U_\alpha$ is also metrizable.
\end{remark}

\begin{proposition}$\label{Thom-Z/2}$ There are Thom isomorphisms in homology
\begin{equation}\label{Thom-isom}H_* (\C_i, \C_{i-1}; \bbZ/2) \isom H_{*-2c(\mu_i)} \left(\C_{\mu_i}; \bbZ/2\right),
\end{equation}
and similarly in the case of a non-orientable surface.
\end{proposition}  

\noindent {\bf Proof.}  In the orientable case, $\cA (E)$ is an (affine) Banach manifold, as are the open subsets $\C_i$.  Daskalopoulos~\cite{Dask} showed that  $\C_{\mu_i}$ is a locally closed submanifold of finite codimension, so by Lemma~\ref{tubular}, $\C_{\mu_i}$ has a tubular neighborhood $\tau_i$ in $\C_i$.
When $\Sigma$ is a non-orientable surface, pulling back to the orientable double cover $\widetilde{\Sigma}$ yields an embedding
$\cA (E)\injects \cA (\widetilde{E})$.  The image is the fixed point set of an involution $\tau$ induced by the deck transformation on $\widetilde{\Sigma}$ (see Ho--Liu \cite{Ho-Liu-non-orient}).   Morse strata in $\cA (\widetilde{E})$ are intersections of Morse strata in $\cA^{k-1} (\widetilde{E})$ with $\cA (E)$, hence are locally closed submanifolds of finite codimension, and we  apply Lemma~\ref{tubular} to obtain a tubular neighborhood $\tau_i$.  
The isomorphisms (\ref{Thom-isom}) come from excising the complement of $\tau_i$ in $\C_i$ and applying the Thom Isomorphism Theorem to the bundle $\tau_i\to \C_{\mu_i}$.
$\hfill \Box$  \vspace{.15in}

%%%%%%%%%%%%%%%%%%%%%%%%%%%%%%%%%%%%%%%%%%%%%%%%%%%
%%%%%%%%%%%%%%%%%%%%%%%%%%%%%%%%%%%%%%%%%%%%%%%%%%%
%%%%%%%%%%%%%%%%%%%%%%%%%%%%%%%%%%%%%%%%%%%%%%%%%%%
%%%%%%%%%%%%%%%%%%%%%%%%%%%%%%%%%%%%%%%%%%%%%%%%%%%

\section{Connectivity of the space of central Yang--Mills connections}$\label{conn}$

Recall that on a bundle $E$ over a Riemann surface $M^g$, the central Yang--Mills connections form the minimum critical set $\Cmin(E)$ of the Yang--Mills functional $L: \cA(E)\to \R$.  The stable manifold of this critical set is the set $\C_{\mathrm{ss}} (E)$ of semi-stable holomorphic structures on $E$, which we refer to as the central stratum, and the Yang--Mills flow provides a deformation retraction $\C_{\mathrm{ss}} (E)\heq \Cmin(E)$.

Using the existence of tubular neighborhoods for Yang--Mills strata, we give a precise formula (Theorem~\ref{conn-calc}) for the connectivity of the spaces $\C_{\mathrm{ss}} (E)\heq \Cmin(E)$, depending only on the genus of $M$ and the rank and Chern number of $E$.  In most cases, we obtain a similar result (Theorem~\ref{conn-calc-non-orient}) for the space $\flatc(E)$ of flat connections on a bundle $E$ over non-orientable surfaces $\Sigma$.  In this case, the connectivity depends only on the genus $\wt{g}$ of the orientable double cover $\wt{\Sigma}$ and the rank of $E$.
Upon considering the homotopy orbits of these spaces under the actions of the gauge groups, 
these results lead to precise formulas for the connectivities of the natural maps from these homotopy orbit spaces to the classifying spaces of the gauge groups (Corollary~\ref{conn-cor}).

The starting point for these calculations is a homological argument, which reduces the problem to a combinatorial question about the codimensions of the Yang--Mills strata.  It is worth noting that the partial ordering $\leqs$ on our strata does \emph{not} respect codimensions (see Example~\ref{non-minimal}); this complicates the argument somewhat.

\begin{proposition}$\label{conn-prop}$ Let $M^g$ be a Riemann surface of genus $g>0$, and let $E$ be a vector bundle over $M$ of rank $n$ and Chern number $k$.  Let $d = d(E)$ denote the minimum (non-zero) codimension of a Harder--Narasimhan stratum in the space $\C(E)$ of holomorphic structures on $E$.  Then the space  $\mathcal{N}_{ss} (E)$ of central Yang--Mills connections on $E$ is precisely $(d-2)$--connected.  

Similarly, let $\Sigma$ be a non-orientable surface and let $E$ be a complex bundle over $\Sigma$.
Let $d=d(E)$ denote the minimum positive codimension of a stratum in the space of connections 
$\cA(E)$.  If $\cA(E)$ contains no strata of codimension $d+1$, then $\flatc(E)$ is precisely $(d-2)$--connected. 
\end{proposition}
 
\noindent {\bf Proof.}  The proofs in the orientable and non-orientable case are essentially identical, so we work in the orientable case (the extra hypothesis in the non-orientable case is automatically satisfied in the orientable case because there the codimensions are always even).

We begin by recalling that by the work of Daskalopoulos~\cite{Dask} and R{\aa}de~\cite{Rade}, the Yang--Mills flow provides a deformation retraction from the space $\Css (E)$ of semi-stable holomorphic structures on $E$ to its critical set $\mathcal{N}_{ss} (E)$, and hence we may work with $\Css (E)$.
Using transversality arguments, it was shown in~\cite[Section 4]{Ramras-surface} that $\pi_i \Css (E) = 0$ for $i \leqs d-2$ (that argument was stated only for the case $k=0$, but works without change in the general case).  We must show that $\pi_{d-1} \Css (E)$ is non-zero.  Since $\Css(E)$ is (at least) $(d-2)$--connected, it suffices, by the Hurewicz Theorem, to prove that $H_{d-1} (\Css(E); \bbZ)\neq 0$.  In fact, we claim that it is enough to show that $H_{d-1} (\Css(E); \bbZ/2\bbZ) \neq 0$.
The Hurewicz Theorem implies that $H_i (\Css(E); \bbZ) = 0$ for $i<d-1$, and hence Tor$(H_{d-2} (\Css (E); \bbZ), \bbZ/2\bbZ) = 0$.  By the Universal Coefficient Theorem, we now have $H_{d-1} (\Css(E); \bbZ/2\bbZ) \isom H_{d-1} (\Css(E); \bbZ) \otimes \bbZ/2\bbZ$, so if $H_{d-1} (\Css(E); \bbZ/2\bbZ)$ is non-zero, we must have $H_{d-1} (\Css(E); \bbZ)\neq 0$ as well.
From now on, all homology groups will be taken with $\bbZ/2\bbZ$--coefficients, and we will drop the coefficient group from our notation.  
\footnote{Our reason for working mod 2 is that in the non-orientable case, the normal bundles to the Yang--Mills strata are \emph{real} vector bundles.  In genus at least 2, these bundles are in fact orientable by Ho--Liu--Ramras~\cite{Ho-Liu-Ramras}, but in genus 1 this is not known.
}  

Let $\prec$ denote the linear ordering on the set of Yang--Mills strata guaranteed by Proposition~\ref{total-order}; we will denote the strata by $\Css = \C_{\mu_0} \prec \C_{\mu_1} \prec \cdots$.
Let 
 $\C_{\mu_m}$ be the first stratum with codimension $d$.  As before, we use the notation
$$\C_j = \bigcup_{i=1}^j \C_{\mu_i}.$$
We claim that $H_{d-1} (\Css) \isom H_{d-1} (\C_{m-1})$.   
By Proposition~\ref{Thom-Z/2}, we have isomorphisms
$$H_*(\C_j, \C_{j-1}) \isom H_{*-\mathrm{codim} (\C_{\mu_j})} (\C_{\mu_j}), \,\,\, j = 1, 2, \ldots, m-1$$
By construction, $\mathrm{codim} (\C_{\mu_j}) > d$ and hence these relative terms are zero in dimensions $d-1$ and $d$.  Thus the long exact sequences of the pairs $(\C_j, \C_{j-1})$ provide  isomorphisms
$$H_{d-1} (\Css) = H_{d-1} (\C_0) \isom H_{d-1} (\C_1) \isom \cdots \isom H_{d-1} (\C_{m-1}).$$

It will now suffice to show that
$H_{d-1} (\C_{m-1}) \neq 0$.  We will argue by contradiction.  
Note that by Lemma~\ref{small-codim}, there are finitely many strata of codimension $d$, say $\C_{\mu_{m_0}}, \ldots, \C_{\mu_{m_l}}$ (with $m=m_0$), and all other strata have codimension at least $d+2$.

Now assume $H_{d-1} (\C_{m-1}) = 0$.  We will prove that $H_{d-1} (\C_{m_i - 1}) = 0$ for $i = 0, 1, \ldots l$.  The proof is by induction on $i$; the base case is our initial assumption.  Now, assuming
 $H_{d-1} (\C_{m_i- 1}) = 0$, consider the long exact sequence in homology for the pair
$(\C_{m_i}, \C_{m_i -1})$.  This sequence has the form
$$0 = H_{d-1} \C_{m_i - 1} \maps H_{d-1} \C_{m_i} \maps H_{d-1} (\C_{m_i}, \C_{m_i - 1})\maps \cdots,$$
and by Proposition~\ref{Thom-Z/2} the last term is zero.  Hence the middle term is zero as well.
Now, since the strata between 
 $\C_{\mu_{m_i}}$ and $\C_{\mu_{m_{i+1}}}$ all have codimension greater than $d$, applying Proposition~\ref{Thom-Z/2} again gives isomorphisms
$$0 = H_{d-1} \C_{m_i} \isom H_{d-1} \C_{m_i + 1} \isom \cdots \isom H_{d-1} \C_{m_{i+1} - 1},$$
 completing the induction.  So we conclude that $H_{d-1} (\C_{m_l-1}) = 0$.  
 
The long-exact sequence for the pair $\C_{m_l - 1} \subset \C_{m_l}$ has the form
$$\cdots \maps H_d (\C_{m_l}) \maps H_d (\C_{m_l}, \C_{m_l - 1}) \maps H_{d-1} \C_{m_l - 1} = 0.$$ 
Since $\C_{\mu_{m_l}}$ has codimension $d$, Proposition~\ref{Thom-Z/2} implies that the relative term is non-zero.  Hence the left-hand term
$H_d (\C_{m_l})$ must be non-zero as well.  But all the remaining strata have codimension at least $d+2$, meaning that  Proposition~\ref{Thom-Z/2}  and Corollary~\ref{gauge-hom}  give isomorphisms
$$H_d (\C_{m_l}) \isom  H_d (\C_{m_l + 1}) \isom \ldots \isom H_d (\C(E)).$$
(Note that here we are only using Corollary~\ref{gauge-hom} with $\bbZ/2\bbZ$--coefficients, in which case the result follows immediately from Proposition~\ref{total-order}.)
Since $\C(E)$ is contractible, this is a contradiction and the proof is complete.
$\hfill \Box$  \vspace{.15in}

The additional hypothesis in the non-orientable case is satisfied in almost all cases, as we will see.
Before beginning the computation of $d(E)$, we note an immediate corollary.  Recall that a map $X\to Y$ is $n$--connected if it induces an isomorphism on $\pi_k$ for $k\leqs n$ and a surjection on $\pi_{n+1}$.
For a rank $n$ bundle $E$ over a surface $M$, the spaces $\Map^E(M, BU(n))$ (the space of classifying maps for $E$) 
and $\Map_*^E(M, BU(n))$ (the subspace of based maps) are models for $B\G(E)$ and $B\G_0 (E)$ respectively~\cite[Section 2]{A-B}.  Hence we obtain fibration sequences 
\begin{equation*}\begin{split} & \Cmin (E) \maps \Cmin (E)_{h\G(E)}\stackrel{q}{\maps} \Map^E(M, BU(n)),\\
\mathrm{and\,\,\,} & \Cmin (E) \maps \Cmin (E)_{h\G_0(E)}\stackrel{q_0}{\maps} \Map_*^E(M, BU(n)),
\end{split}
\end{equation*}
and since the quotient map $\Cmin(E)\to \Cmin(E)/\G_0 (E)$ is a principal bundle~\cite{Mitter-Viallet}, we have weak equivalences 
$$\Cmin (E)_{h\G_0(E)} \heq \Cmin(E)/\G_0 (E) \textrm{  and  }
\Cmin(E)_{h\G} \heq (\Cmin(E)/\G_0 (E))_{hU(n)}.$$  
When $\Cmin(E)$ consists of flat connections, the quotient $\Cmin(E)/\G_0(E)$ is simply the representation space $\Hom(\pi_1 M, U(n))$.  (For details on these issues, we refer the reader to~\cite[Sections 3 and 5]{Ramras-surface}.

By examining the long exact sequences of these fibrations, one sees that since $\pi_{d-1} \Cmin (E) \neq 0$, the maps $q$ and $q_0$ cannot induce isomorphisms on $\pi_{d-1}$ \emph{and} surjections on $\pi_{d}$ (although it is unclear which of these fails).  Thus we have:

\begin{corollary}$\label{conn-cor}$
For any complex vector bundle $E$ over a surface $M$, the maps 
$$\Cmin (E)_{h\G(n)}\to \Map^E(M, BU(n)) \mathrm{\,\,\,\, and \,\,\,\,} \Cmin (E)_{h\G_0(n)} \to \Map_*^E(M, BU(n))$$
and precisely $(d(E) - 2)$--connected, where $d(E)$ is the connectivity of $\Cmin(E)$ and is computed (in nearly all cases) in Theorems~\ref{conn-calc} and~\ref{conn-calc-non-orient}.
\end{corollary}

In the orientable case, the integral (co)homology of $B\G(E)$ and $B\G_0(E)$ were computed by Atiyah and Bott and found to be torsion--free~\cite[Section 2]{A-B}.  For non-orientable surfaces, the rational (co)homology may be computed by the same method.  Hence Corollary~\ref{conn-cor} yields computations of the equivariant (co)homology groups $H^*_{\G(E)} (\Cmin (E))$ and $H^*_{\G_0(E)} (\Cmin(E)) = H^*(\Cmin(E)/\G_0(E))$ below dimension $d(E)-2$.

We now turn to the question of computing the minimum codimension of a non-central stratum.  The case of a trivial bundle was addressed in Ramras~\cite[Lemma 4.5]{Ramras-surface}.  To handle non-trivial bundles, we will need some definitions and lemmas regarding the codimension of the Harder--Narasimhan strata.  This approach will provide an alternate proof in the case $k=0$.

\begin{definition}$\label{codim}$ Let $\mu = ((n_1, k_1), \ldots, (n_r, k_r)) \in \I(n,k)$ be an admissible sequence.  We define 
$$c_1 (\mu) = \sum_{i>j} n_i k_j - n_j k_i \,\,\, \mathrm{and}\,\,\, c_2 (\mu) = (g-1) \sum_{i>j} n_i n_j.$$
Note that the complex codimension of the corresponding Harder--Narasimhan stratum is given by 
$c(\mu) = c_1 (\mu) + c_2 (\mu)$.
\end{definition}

We need some lemmas.  The first follows immediately from the definitions.

\begin{lemma}$\label{r=2}$ For any $n\in \bbN$, $k\in \bbZ$, and $\mu =  ((n_1, k_1),Ê\ldots, (n_r, k_r)) \in \I(n,k)$ with $r>2$, we have
$$c ((n_1, k_1),Ê\ldots, (n_r, k_r)) 
> c \left( \left(\sum_{i=1}^{r-1} n_i, \sum_{i=1}^{r-1} k_i\right), (n_r, k_r)\right).$$
In particular, any admissible sequence minimizing the function $c$ must be of length 2.
\end{lemma}

For $a\in \bbR$, we let $\lceil a \rceil$ will denote the smallest integer \emph{strictly greater} than $a$ (so for $a\in \bbZ$, we set $\lceil a \rceil = a+1$).  This convention will simplify our notation.

\begin{definition}$\label{small-cd-strata}$ For any $n\in \bbN$, $k\in \bbZ$, and $m=1, \ldots, n-1$,
let
$$\mu_{m} = \mu_m (n,k) = ((m, \lceil \frac{k m}{n} \rceil), (n-m, k - \lceil \frac{km}{n} \rceil)).$$
We define
$$\I'(n,k) = \{\mu_m \,\, : \,\, 0<m<n\} \subset \I(n,k).$$
(The line from $(0,0)$ to $(n, k)$ passes through $(m, \frac{km}{n})$, so $\mu_m$ is admissible.)
\end{definition}

\begin{lemma}$\label{ceiling}$
For any $n\in \bbN$, $k\in \bbZ$, and any admissible sequence 
$$\mu =  ((m, l),Ê(n-m, k- l)) \in \I(n,k)$$
of length two, we have $c (\mu) \geqs c (\mu_m),$
with equality holding only when $\mu = \mu_m$.

Hence if $\mu\in \I(n,k)$ minimizes the function $c$, then $\mu \in \I'(n,k)$.
\end{lemma}

We now consider what values the function $c_1$ may take on the set $\I'(n,k)$.  

\begin{definition}  Given an integer $r$ and a natural number $n$, we let $[r]_n$ denote the unique integer between $1$ and $n$ satisfying $r\equiv [r]_n$ (mod $n$).
\end{definition}

\begin{lemma}$\label{cyclic}$
For any $n\in \bbN$ and $k\in \bbZ$, we have $c_1 (\mu_m) = [km]_n$, and  
$$c_1 (\I'(n,k)) = \{\mathrm{gcd}(n,k), 2\mathrm{gcd}(n,k), \ldots, n\}$$
(unless $\mathrm{gcd}(n,k) = 1$, in which case $n$ is not included in this set).
\end{lemma}
\noindent {\bf Proof.}  Let $\lceil \frac{km}{n} \rceil = \frac{km}{n} + \frac{\epsilon_m}{n}$, 
and note that $\epsilon_m = n\lceil \frac{km}{n} \rceil - km \equiv -km$ (mod $n$).  Since 
$\epsilon_m$ is an integer between $1$ and $n$, we have $\epsilon_m = [-km]_n$.

Now, for any sequence
$\mu_m \in \I'(n,k)$, we have 
\begin{equation*}
\begin{split}
c_1 (\mu_m) 
&= c_1 ((m, \frac{km}{n}+ \frac{\epsilon_m}{n}), (n-m, k- \frac{km}{n} - \frac{\epsilon_m}{n}))\\
& = (n-m)(\frac{km}{n}+ \frac{\epsilon_m}{n}) - m (k- \frac{km}{n} - \frac{\epsilon_m}{n})\\
& = \epsilon_m = [-km]_n,
\end{split}
\end{equation*}
as desired.
Now, consider the set $\{1, \ldots, n\}$ as a cyclic group under addition modulo $n$.  
Then $c_1 (\I'(n,k)) = \{[-km]_n : m = 1, 2, \ldots, n-1\}$.  Since $[-mk]_n \equiv m[-k]_n$ (mod $n$), $c_1 (\I'(n,k))$ is the subgroup of $\{1, \ldots, n\}$ generated by $[-k]_n$.  But $\gcd(n,k)$ is the minimal element of this subgroup, so the lemma follows.
$\hfill \Box$  \vspace{.15in}

We can now determine the exact connectivity of the space of central Yang--Mills connections on any bundle $E$ over a Riemann surface.  Recall that by Proposition~\ref{conn-prop}, it suffices to calculate the minimum codimension of a non-central Harder--Narasimhan stratum in the space $\C (E)$.

\begin{theorem}$\label{conn-calc}$
If $E$ is a complex vector bundle over a Riemann surface $M^g$ ($g>0$), then the connectivity of the space $\Cmin (E)$ of central Yang--Mills connections is given by
$$d(E) - 2=\left\{ \begin{array}{ll}
					2\gcd (n,k) - 2 \,\,\,\,\,\,\,\,\,\,\,\,\,\,\,\, &\mathrm{if} \, g=1\\
					2\min ([k]_n, [-k]_n) + 2(g-1)(n-1) - 2& \mathrm{if} \, g>1
				   \end{array}		
			\right.
$$  
\end{theorem}
\noindent {\bf Proof.} We know (Lemma~\ref{ceiling}) that $c$ is minimized by $c(\mu_m)$ for some integer $m$ between $1$ and $n-1$.   For $g=1$, the function $c_2$ vanishes, so the result follows from Lemma~\ref{cyclic}.
For $g>1$, we begin by noting that
\begin{equation}\label{cm-c1}
\begin{split}
c(\mu_m) - c(\mu_1) &= c_1(\mu_m) - c_1 (\mu_1) + c_2(\mu_m) - c_2 (\mu_1) \\
				&= [-mk]_n - [-k]_n + (g-1)m(n-m) - (g-1)(n-1)\\
				&= [-mk]_n - [-k]_n + (g-1)(mn-m^2 -n+1).
\end{split}
\end{equation}

Assuming $n\geqs 6$, we will show that if $m\leqs n/2$
then $c(\mu_m)\geqs c(\mu_1)$, and if $m\geqs n/2$ then $c(\mu_m)\geqs c(\mu_{n-1})$.  This will suffice to prove the theorem for $n\geqs 6$, since $c(\mu_1) = [-k]_n + (g-1)(n-1)$ and $c(\mu_{n-1}) = [k]_n + (g-1)(n-1)$.  The cases $n<6$ can be checked by hand.

We first consider the case $2\leqs m \leqs n/2$; we may assume $n>2$.
The function $f_n(m) = mn-m^2-n+1$ has derivative $(f_n)^{'}(m) = n-2m\geqs 0$ (since $m\leqs n/2$)
and hence this function is minimized at $m=2$, where we have 
$f_n(2) = 2n-4-n+1 = n-3\geqs 0$.  So $mn-m^2 -n+1$ is always positive, and hence
\begin{equation}\label{c2}
(g-1)(mn-m^2 -n+1) \geqs mn-m^2-n+1.
\end{equation}
Equations (\ref{cm-c1}) and (\ref{c2}) imply that
\begin{equation}\label{cmc1c2}
c(\mu_m) - c(\mu_1) \geqs [-mk]_n - [-k]_n + mn-m^2-n+1.
\end{equation}

For later reference, we work with a generic integer $r$ in place of $-k$.  Note that if $\frac{l}{m} n< [r]_n\leqs \frac{l+1}{m} n$ ($l = 0, 1, \ldots, m-1$), then $0<m[r]_n - l n \leqs n$, so 
we have
\begin{equation}\label{r}
[mr]_n = m[r]_n - ln.
\end{equation}
Furthermore, since $l\leqs m-1$ and $m\leqs n/2$, (\ref{r}) implies 
\begin{equation}\label{l}
\begin{split}
[mr]_n - [r]_n &= (m[r]_n - ln) - [r]_n = (m-1)[r]_n - ln \\
			    & > (m-1)\frac{ln}{m} - ln = -\frac{ln}{m}\\ 
			    & \geqs -\frac{(m-1)n}{m} = -n + n/m\\
			    & \geqs -n+\frac{n}{n/2} = -n+2.
\end{split}
\end{equation} 
Combining (\ref{l}) and (\ref{cmc1c2}) yields
\begin{equation}\label{g(m)}
\begin{split}
c(\mu_m) - c(\mu_1) &> -n+2 + mn-m^2-n+1\\
				&=n(m-2) -m^2 + 3.
\end{split}
\end{equation}
Letting $h_n(m) = n(m-2) -m^2 + 3$, we have $(h_n)'(m) = n-2m>0$ (since $m\leqs n/2$).
On the interval $3\leqs m\leqs n/2$, $h_n (m)$ is minimized at $m=3$, so
$$c(\mu_m) - c(\mu_1) > h_n(3) = n - 9 + 3 \geqs 0$$
since $n\geqs 6$.  Note that $h_n (2) = -1$, so a different estimate is needed when $m=2$.

When $m=2$, we have
$\frac{l}{2} n< [r]_n\leqs \frac{l+1}{2} n$ for either $l=0$ or $l=1$.  By (\ref{r}),
\begin{equation}\label{l2}
\begin{split}
[2r]_n - [r]_n &= (2[r]_n - ln) - [r]_n = [r]_n - ln \\
			    & > \frac{ln}{2} - ln = -\frac{ln}{2} \geqs -\frac{n}{2}
\end{split}
\end{equation}
Combining (\ref{cmc1c2}) and (\ref{l2}) yields (for $n\geqs 6$)
$$c(\mu_2) - c(\mu_1) > -n/2 + 2n - 4 - n + 1 = n/2 - 3\geqs 0.$$

Thus we have shown that $c(\mu_m)\geqs c(\mu_1)$ for $2\leqs m\leqs n/2$.
The proof that $c(\mu_{n-m})\geqs c(\mu_{n-1})$ for $2\leqs m\leqs n/2$ is symmetrical:
as before we have
$$c(\mu_{n-m}) - c(\mu_{n-1}) \geqs [mk]_n - [k]_n + mn-m^2-n+1.$$
When $m=2$ (\ref{l2}) implies that $[2k]_n - [k]_n\geqs -n/2$, and in general (\ref{l}) implies
that $[mk]_n - [k]_n\geqs -n + 2$.  The argument now proceeds identically.
$\hfill \Box$  \vspace{.15in}

When $k=0$, $[0]_n = n$ and hence the formula given here recovers that found in~\cite[Lemma 4.5]{Ramras-surface}.
We also note that when $g=1$ and $k\neq 0$, this result shows that the connectivity of the space of central Yang--Mills connections does not tend to infinity with $n$.  
The following example shows that the strata of minimal codimension are not, in general, minimal covers of the central stratum.  

\begin{example}$\label{non-minimal}$
When $n=6$, $k=2$, and $g>1$, Theorem~\ref{conn-calc} shows that the stratum
$\mu_1 = ((1, 1), (5, 1))$ has minimum (complex) codimension, given in this case by $4+5(g-1)$.  However, this stratum lies above the stratum
$((2, 1), (4, 1))$, which has complex codimension $2 + 8(g-1)$.
\end{example}

Since the critical values of the Yang--Mills functional are given by 
Proposition~\ref{inf} (see Remark~\ref{l(mu)}), one can show by a combinatorial argument that these critical values respect the partial ordering on the strata (and of course one may check this directly in Example~\ref{non-minimal}).  Hence the Yang--Mills functional is not self-indexing, even after scaling.

We now turn to the case of a complex vector bundle $E$ over a non-orientable surface $\Sigma$.  Here the minimum critical set of the Yang--Mills functional is the space $\flatc(E)$ of flat connections.  A combinatorial argument (simpler than the previous one) allows us to calculate the connectivity of $\flatc(E)$ in most cases.

\begin{theorem}$\label{conn-calc-non-orient}$
Let $E$ be a complex bundle of rank $n>1$ over a non-orientable surface $\Sigma$, and let $\wt{g}$ denote the genus of the orientable double cover $\wt{\Sigma}$.  If $\wt{g}\geqs 2$ and $n\geqs 9$, then the space $\flatc(E)$ of flat connections on $E$ is precisely $(2n\wt{g} - 3 \wt{g} - 1)$--connected.
\end{theorem}
\noindent {\bf Proof.} We will show that the minimum (positive) codimension of a Yang--Mills stratum in the space $\A(E)$ of connections on $E$ is precisely two more than the stated connectivity.  Moreover, we will show that any other positive-codimension stratum has codimension at least two more than the minimum; the result then follows from Proposition~\ref{conn-prop}.  We will point out the differences in genus 1.

To begin, recall from Ho and Liu~\cite{Ho-Liu-non-orient} or Ho--Liu--Ramras~\cite{Ho-Liu-Ramras} that $\A(E)$ embeds as the set of fixed points of an involution on $\A(\wt{E})$, and each Yang--Mills stratum in $\A(E)$ is the collection of fixed points lying inside some given Yang--Mills stratum of $\A(\wt{E})$.  In fact, any stratum in $\A(\wt{E})$ containing fixed points corresponds to an admissible sequence of the form
\begin{equation}\label{mu-no}
\mu = ((n_1, k_1), \ldots, (n_r, k_r), (n_0, 0), (n_r, -k_r), \ldots, (n_1, -k_1)),
\end{equation}
where $\sum n_i = n$  and the bundle $\wt{E}$ is necessarily trivial (although not all such strata contain fixed points~\cite[Section 7.1]{Ho-Liu-non-orient}).  We will call such sequences \emph{symmetric}.  The set of fixed points lying inside a symmetric stratum, if non-empty, has real codimension $c(\mu)$ inside $\A(E)$, where $c(\mu)$ is the \emph{complex} codimension of the stratum $\A_\mu$ inside $\A(\wt{E})$ and is given by the formula in Definition~\ref{codim}.

In analogy with Lemma~\ref{r=2}, one sees that with $\mu$ as in (\ref{mu-no}),
$$c(\mu) \geqs c\left( \left(\sum n_i, \sum k_i\right), (n_0, 0), \left(\sum n_i, -\sum k_i\right) \right) + 2$$
when $r>1$,
and hence any symmetric stratum minimizing $c$ must be of the form
$$\mu = ((n_1, k_1), (n_0, 0), (n_1, -k_1)).$$
Next, it is again elementary to check that 
$$c((n_1, k_1), (n_0, 0), (n_1, -k_1)) \geqs c((n_1, 1), (n_0, 0), (n_1, -1)) + 2$$
for $k_1>1$, so the minimum codimension can only be achieved by the strata
$$\mu_i  = ((i, 1), (n_0, 0), (i, -1)),$$
$i = 1, \ldots, \lfloor n/2 \rfloor$ (here $\lfloor n/2 \rfloor = n/2$ if $n$ is even and $\lfloor n/2 \rfloor = (n-1)/2$ if $n$ is odd).
Now $c(i):= c(\mu_i) = 2n - 2i + (2ni - 3i^2)(\wt{g}-1)$ is quadratic in $i$ with a maximum at $\frac{n-1/(\wt{g}-1)}{3}$, and it is elementary to check that for $\wt{g}>1$, 
$$c(2) \geqs c(1) + 2  \mathrm{\,\, for\,\,} n\geqs 9 \mathrm{\,\,\,\,\,\,\,\, and \,\,\,\,\,\,\,\,} c(\lfloor n/2 \rfloor) \geqs  c(n/2) \geqs c(1) + 2 \mathrm{\,\, for\,\,} n\geqs 12.$$
Hence when $\wt{g}>1$ and $n\geqs 12$, the $c(i)$ is minimized when $i=1$ (and $c(1) = 2n\wt{g} - 3 \wt{g} + 1$).  The stratum $\mu_1 = ((1, 1), (n-1, 1), (1,-1)$ is in fact non-empty because the degree zero factor has dimension $n-1> 0$~\cite[Section 7.1]{Ho-Liu-non-orient}.   Hence $c(1)$ gives the minimum 
minimum positive codimension of a non-empty stratum.  The cases $\wt{g}>1$, $9\leqs n\leqs 11$ can be checked by hand.  
$\hfill \Box$  \vspace{.15in}

Our reduction to the strata $\mu_i$ did not require $n\geqs 9$, so the remaining cases may be computed by hand. In most cases, one still obtains the connectivity of $\flatc(E)$ precisely.  But when $n=5$ and $\wt{g} = 2$ or $4$, there is a (non-empty) stratum of codimension one more than the minimum, so Proposition~\ref{conn-prop} does not apply.  

The case $\wt{g}=1$, where $\Sigma = K$ is the Klein bottle, is slightly different.  When $n$ is even, strata of the form $((n/2, k), (n/2, -k))$ in $\A(E)$ may be empty: each such stratum for the trivial bundle on the double cover $S^1\cross S^1$ contains connections from either $E^+$ or $E^-$, but not both.  When $2(n/2) + k + 2 = n+k+2$ is even, this stratum contains connections from $E^+$, and when $n+k+2$ is odd, it contains connections from $E^-$~\cite[Proposition 7.1]{Ho-Liu-non-orient}.

The reductions in the proof of Theorem~\ref{conn-calc-non-orient} show that for odd $n$, the minimum codimension of a non-empty stratum is $n+1$ and the connectivity of $\flatc(E)$ is precisely $n-1$.  If $n$ is even and $E = E^-$ is the non-trivial bundle, one finds that the stratum $((n/2, 1), (n/2, -1))$ is non-empty and again gives the minimum codimension, namely $n$.  So $\flatc(E^-)$ is precisely $(n-2)$--connected.  Finally, if $E^+ = K\cross \bbC^n$, then these reductions show that any non-empty stratum has codimension at least two more than either
 $((n/2, 2), (n/2, -2))$ or $\mu_{n/2 -1} = ((n/2-1, 1), (2, 0), (n/2 -1, -1))$, both of which are non-empty.  The minimum codimension of a non-empty stratum is thus $c(n/2 - 1) = n+2$, and $\flatc(E^+)$ is precisely $n$-connected.

%%%%%%%%%%%%%%%%%%%%%%%%%%%%%%%%%%%%%%%%%%%%%%%%%%%
%%%%%%%%%%%%%%%%%%%%%%%%%%%%%%%%%%%%%%%%%%%%%%%%%%%

\end{document}